\documentclass[twoside,emlines,bezier,12pt]{article}
\usepackage{amsfonts,amssymb}

\begin{document}

\begin{center}

{\Large  Order-indices and order-periods of $3\times 3$ matrices
over commutative inclines}

$ $

\bigskip {Song-Chol Han$\textrm{*}$,\quad Gum-Song Sin}

\vskip 3mm

{\small \centerline{{\it Faculty of Mathematics, Kim Il Sung
University, Pyongyang, DPR Korea}}}

\renewcommand{\thefootnote}{\alph{footnote}}
\setcounter{footnote}{-1} \footnote{$\textrm{*}$ Corresponding
author. {\it E-mail:} ryongnam5@yahoo.com}
\end{center}

\vspace{0.3cm}

{\small \noindent {\bf Abstract}
\medskip

An incline is an additively idempotent semiring in which the
product of two elements is always less than or equal to either
factor. By making use of prime numbers, this paper proves that
$A^{11}\leqslant A^5$ for all $3\times 3$ matrices $A$ over an
arbitrary commutative incline, thus giving an answer to an open
problem ``For $3\times 3$ matrices over any incline (even
noncommutative) is $X^5\geqslant X^{11}$ ?'', proposed by Cao, Kim
and Roush in a monograph {\it Incline Algebra and Applications},
{\it 1984}.
\medskip

\medskip

}
\medskip

$ $

{\section*{\noindent \bf 1. Introduction}}

Inclines are additively idempotent semirings in which products are
less than or equal to either factor. Boolean algebra, fuzzy
algebra and distributive lattice are examples of inclines. The
study of inclines is generally acknowledged to have started by
Z.Q. Cao in a series of his papers in the early 1980's. Incline
algebra and incline matrix theory have been extensively studied by
many authors. Nowadays, one may clearly notice a growing interest
in developing the algebraic theory of inclines and their numerous
significant applications in diverse branches of mathematics and
computer science such as automata theory, graph theory,
informational systems, complex systems modelling, decision-making
theory, dynamical programming, control theory, clustering and so
on (see \cite{Kim04}). Inclines are also called simple semirings
(see \cite{Golan99}).

Cao et al. \cite{Cao84} introduced the notion of the order-index
and order-period of an element in a partially ordered semigroup,
and proposed an open problem ``For $3\times 3$ matrices over any
incline (even noncommutative) is $X^5\geqslant X^{11}$ ?'' (see
the fourth open problem in Section 5.5 in \cite{Cao84}).

Han and Li \cite{Han07} proved that $A^{k+d}\leqslant A^{k}$ and
thus the order-index of $A$ is at most $k$ for all $n\times n$
matrices $A$ over an arbitrary commutative incline, where
$k=(n-1)^2+1$ and $d$ is any given multiple of $[n]$ satisfying
$d\geqslant nk(\sum_{i=1}^n P_n^i)$. In the case of $n=3$, one can
easily conclude that $A^{233}\leqslant A^{5}$ and so the
order-index of $A$ is at most 5 for all $3\times 3$ matrices $A$
over any commutative incline (because $k=5$, $[n]=6$, and
$\sum_{i=1}^n P_n^i=P_3^1+P_3^2+P_3^3=15$).

In this paper, we prove that $A^{11}\leqslant A^5$ for all
$3\times 3$ matrices $A$ over an arbitrary commutative incline by
using prime numbers, thereby giving an answer to the above open
problem.
\medskip

$ $

{\section*{\noindent \bf 2. Preliminaries}}

We recall some known definitions and facts.
\medskip

\noindent {\bf Definition 2.1} \cite{Cao84}.\ \ Let $+$ and
$\cdot$ be two binary operations on a nonempty set $L$. An
algebraic system $(L,+,\cdot)$ is called an incline if it
satisfies the following conditions:
\medskip

\noindent (1)\ \ $(L,+)$ is a semilattice,

\noindent (2)\ \ $(L,\cdot)$ is a semigroup,

\noindent (3)\ \ $x(y+z)=xy+xz$ and $(y+z)x=yx+zx$ for all
$x,y,z\in L$,

\noindent (4)\ \ $x+xy=x+yx=x$ for all $x,y\in L$.
\medskip

On an incline $L$, define a relation $\leqslant$ by $x\leqslant
y\Leftrightarrow x+y=y$. It is easy to see that $\leqslant$ is a
partial order on $L$ and that for any $x,y\in L$, the sum $x+y$ is
the least upper bound of $\{x,y\}\subseteq L$, i.e., $x+y=x\vee y$
in the poset $(L,\leqslant)$. It follows that $xy\leqslant x$ and
$yx\leqslant x$ for all $x,y\in L$ and that for any $x,y,z\in L$,
$y\leqslant z$ implies $xy\leqslant xz$ and $yx\leqslant zx$.

An incline $L$ is said to be commutative if $xy=yx$ for all
$x,y\in L$.

The Boolean algebra $(\{0,1\},\vee,\wedge)$ is an incline. More
generally, every distributive lattice is an incline. Each fuzzy
algebra $([0,1],\vee,T)$ is an incline, where $T$ is a t-norm. The
tropical algebra $(\mathbb{R}^{+}_0\cup\{\infty\},\wedge,+)$ is an
incline, where $\mathbb{R}^{+}_0$ is the set of all nonnegative
real numbers.

In the sequel, $L$ always denotes any given commutative incline,
$n$ denotes any given positive integer greater than or equal to
$2$, $\underline{n}$ stands for the set $\{1,2,\ldots,n\}$, and
$[n]$ denotes the least common multiple of integers
$1,2,\ldots,n$. For a nonnegative integer $l$, $\underline{l}^0$
denotes the set of integers $0$ through $l$.

We denote by $L^{n\times n}$ the set of all $n\times n$ matrices
over $L$. Given $A=(a_{ij})\in L^{n\times n}$ and $B=(b_{ij})\in
L^{n\times n}$, we define the product $A\cdot B\in L^{n\times n}$
by $A\cdot B:=\Big(\sum_{k\in \underline{n}}a_{ik}b_{kj}\Big)$.
And we denote $A\leqslant B$ when $a_{ij}\leqslant b_{ij}$ for all
$i,j\in \underline{n}$.

Then $(L^{n\times n},\leqslant,\cdot)$ is a partially ordered
semigroup, i.e., for all $A,B,C,D\in L^{n\times n}$, the following
hold.
\medskip

\noindent (1)\ \ $(AB)C=A(BC)$,

\noindent (2)\ \ $A\leqslant B$ and $C\leqslant D$ $\Rightarrow$
$AC\leqslant BD$.
\medskip

Given a matrix $A\in L^{n\times n}$, its powers are defined as
follows: $A^1:=A,\ A^l:=A^{l-1}A$ for all $l\geqslant 2.$
\medskip

\noindent {\bf Definition 2.2} \cite{Cao84}.\ \ Let $S$ be a
partially ordered semigroup and $a\in S$. The order-index of $a$
is the least such positive integer $k$ that $a^{k+d}\leqslant a^k$
for some positive integer $d$. The order-period of $a$ is the
least such positive integer $d$ that $a^{k+d}\leqslant a^k$ for
some positive integer $k$.
\medskip

\noindent {\bf Definition 2.3} \cite{Han07}.\ \ Let
$V:v_0,v_1,\ldots,v_l$ be a sequence of positive integers such
that $v_i\in \underline{n}$ for all $i\in \underline{l}^0$. We
call $V$ a walk on $\underline{n}$, $l(V):=l$ the length of $V$,
and $v_i$ $(i\in \underline{l}^0)$ the terms of $V$. When
$l\geqslant 1$, for any $p,q\in \underline{n}$, we put
$m(V;p,q):=|\{i\in \underline{(l-1)}^0 \mid v_i=p, v_{i+1}=q\}|$.
Let $U:u_0,u_1,\ldots,u_h$ be another walk on $\underline{n}$. $U$
is called a reduction of $V$ if $u_0=v_0$, $u_h=v_l$ and
$m(U;p,q)\leqslant m(V;p,q)$ for all $p,q\in \underline{n}$.
\medskip

$ $

{\section*{\noindent \bf 3. Main results}}

\noindent {\bf Lemma 3.1.}\ \ {\it Let $P=(p_{ij})$ be an $n\times
n$ matrix consisting of $n^2$ distinct prime numbers. Let
$V:v_0,v_1,\ldots,v_l$ be a walk on $\underline{n}$ with length
$l\geqslant 1$ and let $U:u_0,u_1,\ldots,u_h$ be another walk on
$\underline{n}$ with length $h\geqslant 1$ satisfying $u_0=v_0$
and $u_h=v_l$. If the product $\prod_{r=0}^{h-1} p_{u_{r}u_{r+1}}=
p_{u_{0}u_{1}}p_{u_{1}u_{2}}\cdots p_{u_{h-1}u_{h}}$ divides the
product $\prod_{s=0}^{l-1} p_{v_{s}v_{s+1}}=
p_{v_{0}v_{1}}p_{v_{1}v_{2}}\cdots p_{v_{l-1}v_{l}}$, then $U$ is
a reduction of $V$}.
\medskip

\noindent {\bf Proof.}\ \ Since the entries $p_{ij}$ of $P$ are
distinct primes, the divisibility implies that the multiplicity of
each prime factor $p_{u_{r}u_{r+1}}$ in the product
$\prod_{r=0}^{h-1} p_{u_{r}u_{r+1}}$ is less than or equal to its
multiplicity in the product $\prod_{s=0}^{l-1} p_{v_{s}v_{s+1}}$.
Hence we have $m(U;u_{r},u_{r+1})\leqslant m(V;u_{r},u_{r+1})$ for
all $r\in \underline{(h-1)}^0$, and so $U$ is a reduction of $V$.\
\ $\square$
\medskip

\noindent {\bf Lemma 3.2.}\ \ {\it Every walk on $\underline{3}$
with length $11$ has a reduction with length $5$}.
\medskip

\noindent {\bf Proof.}\ \ Putting
$$
P=\left(
\begin{array}{ccc}
2&3&5\\
7&11&13\\
17&19&23\\
\end{array}
\right)
$$
and using Lemma 3.1, the statement is verified by a direct
computation.\ \ $\square$
\medskip

\noindent {\bf Theorem 3.3.}\ \ {\it If $A\in L^{3\times 3}$, then
$A^{11}\leqslant A^{5}$}.
\medskip

\noindent {\bf Proof.}\ \ Let $A=(a_{ij})$. We denote
$A^{5}=\Big(a^{(5)}_{ij}\Big)$ and
$A^{11}=\Big(a^{(11)}_{ij}\Big)$. For any $i,j\in \underline{3}$,
we have $a_{ij}^{(11)}=\sum_{v_1,v_2,\ldots,v_{10}\in
\underline{3}} a_{iv_1}a_{v_1v_2}\cdots a_{v_{10}j}$. Consider
every summand $a_{iv_1}a_{v_1v_2}\cdots a_{v_{10}j}$ of
$a_{ij}^{(11)}$. By Lemma 3.2, the walk
$i,v_1,v_2,\ldots,v_{10},j$ on $\underline{3}$ with length $11$
has a reduction $i,u_1,u_2,u_{3},u_{4},j$ with length $5$.
Noticing that $L$ is a commutative incline, we obtain
$$
a_{iv_1}a_{v_1v_2}\cdots a_{v_{10}j}\leqslant
a_{iu_1}a_{u_1u_2}\cdots a_{u_{4}j}\leqslant
\sum_{k_1,k_2,k_3,k_{4}\in \underline{3}} a_{ik_1}a_{k_1k_2}\cdots
a_{k_{4}j}=a_{ij}^{(5)}.
$$
Since this inequality holds for all summands of $a_{ij}^{(11)}$,
we have $a_{ij}^{(11)}\leqslant a_{ij}^{(5)}$, as required.\ \
$\square$
\medskip

\noindent {\bf Corollary 3.4.}\ \ The order-periods of $3\times 3$
matrices over any commutative incline are at most 6.
\medskip

$ $

{\section*{\noindent \bf 4. An algorithm for verification}}

We present an algorithm for verifying Lemma 3.2:
\medskip

\noindent {\it Step} 1. Input the $3\times 3$ matrix $P=(p_{ij})$
given in the proof of Lemma 3.2 and put $cnt=0$.

\noindent {\it Step} 2. Choose a walk
$i_0,i_1,\ldots,i_{10},i_{11}$ on $\underline{3}$ with length 11
and compute the product $X=p_{i_0i_1}*p_{i_1i_2}*\cdots
*p_{i_{10}i_{11}}$.

\noindent {\it Step} 3. Choose a walk $i_0,j_1,\ldots,j_4,i_{11}$
on $\underline{3}$ with length 5 and compute the product $Y=p_{i_0
j_1}*p_{j_1 j_2}*\cdots *p_{j_4 i_{11}}$.

\noindent {\it Step} 4. Check if $Y$ divides $X$. If yes, then
output the walk $i_0,i_1,\ldots,i_{10},i_{11}$ and its reduction
$i_0,j_1,\ldots,j_4,i_{11}$, and go to Step 7. Otherwise, go to
Step 5.

\noindent {\it Step} 5. Check if all such walks on $\underline{3}$
with length 5 have been tested. If yes, then go to Step 6.
Otherwise, go to Step 3.

\noindent {\it Step} 6. Output the walk
$i_0,i_1,\ldots,i_{10},i_{11}$ which has no reduction with length
5 and add 1 to $cnt$.

\noindent {\it Step} 7. Check if all walks on $\underline{3}$ with
length 11 have been tested. If yes, then go to Step 8. Otherwise,
go to Step 2.

\noindent {\it Step} 8. Output the value of $cnt$, i.e., the total
number of such walks on $\underline{3}$ with length 11 that have
no reduction with length 5. Output the run time.
\medskip

We describe a computer program corresponding to the algorithm
presented above:
\medskip

{\tt

\#include $\langle$stdio.h$\rangle$

\#include $\langle$time.h$\rangle$

typedef \_\_int64 INT;
\medskip

INT X, Y;

int p[3][3] = \{2, 3, 5, 7, 11, 13, 17, 19, 23\};

int cnt, walk[12], walk\_[6];
\medskip

int g(int k, int m, int last)

\{

    walk\_[k] = m;

    if (k == 4)

    \{

        walk\_[5] = last;

        Y *= p[m][last];

        int ok = false;

        if (X \% Y == 0)

        \{

            ok = true;

            int j;

            for (j = 0; j < 12; j++) printf("\%d ", walk[j] + 1);

            for (printf(":"), j = 0; j < 6; j++) printf(" \%d", walk\_[j] + 1);

            puts("");

        \}

        Y /= p[m][last];

        return ok;

    \}

    for (int i = 0; i < 3; i++)

    \{

        Y *= p[m][i];

        int ok = g(k + 1, i, last);

        Y /= p[m][i];

        if (ok) return 1;

    \}

    return 0;

\}
\medskip

int f(int k, int m)

\{

    walk[k] = m;

    if (k == 11)

    \{

        if (!g(0, 0, m))

        \{

            cnt++;

            for (int j = 0; j < 12; j++) printf("\%d ", walk[j] + 1);

            puts(": no reduction of length 5");

        \}

        return 0;

    \}

    for (int i = 0; i < 3; i++)

    \{

        X *= p[m][i];

        f(k + 1, i);

        X /= p[m][i];

    \}

    return 0;

\}
\medskip

int main()

\{

    freopen("out.txt", "w", stdout);

    double \_t = clock();

    X = Y = 1;

    f(0, 0);

    printf("count = \%d$\backslash$n", cnt);

    printf("\%lf$\backslash$n", clock() - \_t);

    return 0;

\}

}
\medskip

The program has been run in 2 seconds on Windows XP version 2002,
IBM-PC Computer with Intel(R) Pentium(R) 4 CPU 2.80 GHz, 248MB RAM
and 80GB Hard Disk. The displayed result showed that the total
number of such walks on $\underline{3}$ with length 11 that have
no reduction with length 5 is {\it zero}.
\medskip

$ $

{\section*{\noindent \bf 5. Conclusions}}

This paper proved that $A^{11}\leqslant A^5$ for all $3\times 3$
matrices $A$ over any commutative incline. The following problem
is still open: For $3\times 3$ matrices $X$ over any
noncommutative incline is $X^5\geqslant X^{11}$ ?
\medskip

\vspace{0.8cm}


{\renewcommand\baselinestretch{1}
\renewcommand\refname}

\end{document}